\documentclass[12pt]{article}
\usepackage{amssymb,amsmath,mathrsfs}
\usepackage{hyperref}
\usepackage{graphicx}
\usepackage{color}
\usepackage[OT2,OT1]{fontenc}

\def\Sha{\text{\fontencoding{OT2}\selectfont SH}}

\begin{document}

\title{\LARGE\bf Representation of the Fourier transform as a weighted sum of the complex error functions}

\author{
\normalsize\bf S. M. Abrarov\footnote{\scriptsize{Dept. Earth and Space Science and Engineering, York University, Toronto, Canada, M3J 1P3.}}\, and B. M. Quine$^{*}$\footnote{\scriptsize{Dept. Physics and Astronomy, York University, Toronto, Canada, M3J 1P3.}}}

\date{August 5, 2015}
\maketitle

\begin{abstract}
In this paper we show that a methodology based on a sampling with the Gaussian function of kind $h\,{e^{ - {{\left( {t/c} \right)}^2}}}/\left( {{c}\sqrt \pi  } \right)$, where ${c}$ and $h$ are some constants, leads to the Fourier transform that can be represented as a weighted sum of the complex error functions. Due to remarkable property of the complex error function, the Fourier transform based on the weighted sum can be significantly simplified and expressed in terms of a damping harmonic series. In contrast to the conventional discrete Fourier transform, this methodology results in a non-periodic wavelet approximation. Consequently, the proposed approach may be useful and convenient in algorithmic implementation.
\vspace{0.25cm}
\vspace{0.25cm}
\end{abstract}

\section{Introduction}

Sampling is an efficient methodology that can be used for functional approximations \cite{Bracewell2000, Hansen2014}. It is based on interpolation through a set of discrete points $\left\{ {{t_n},f\left( {{t_n}} \right)} \right\}$ selected from a function $f\left( t \right)$ within an interval $\left[ {a,b} \right]$. The sinc function
$$
{\rm{sinc}}\left( t \right) = \left\{ 
\begin{aligned}
&\frac{{\sin t}}{t}, \quad\,\, t \ne 0\\
&1, \qquad\quad t = 0,
\end{aligned} \right.
$$
is one of the most popular in sampling \cite{Stenger2011, Rybicki1989, Lether1998}. Specifically, due to its remarkable sampling property $f\left( t \right)$ can be represented as a series approximation \cite{Rybicki1989, Lether1998}
\begin{equation}\label{eq_1}
f\left( t \right) \approx \sum\limits_{n =  - N}^N {{\rm{sinc}}\left( \frac{\pi}{h} \left( {t - {t_n}} \right ) \right)f\left( {{t_n}} \right)},
\end{equation}
where $h$ is a small parameter that at $t_{n} = n h$ can be regarded as the step. We have shown recently a new methodology of sampling based on incomplete cosine expansion of the sinc function in numerical integration \cite{Abrarov2015a}.

Another very useful sampling tool can be constructed by using the Gaussian function ${e^{ - {t^2}}}$. The sampling property of the Gaussian function ${e^{ - {t^2}}}$ can be observed by changing its variable as $t \to t/\varepsilon $ and then considering the limit 
$$
\mathop {\lim }\limits_{\varepsilon  \to 0} {e^{ - {{\left( {t/\varepsilon } \right)}^2}}} = \delta \left( t \right),
$$
where the Kronecker\text{'}s delta $\delta \left( t \right)$ is a building block of {\it{the shah function}}, also known as {\it{the comb function}}:
$$
\Sha\left( t \right) = \sum\limits_{n =  - \infty }^\infty  {\delta \left( {t - {t_n}} \right)}
$$
(see application of {\it{the shah function}} $\Sha\left( t \right)$ in sampling in \cite{Bracewell2000}). Consequently, this example shows that the Gaussian function can also be utilized for functional approximations by sampling methodology.

It is convenient to modify the Gaussian function ${e^{ - {t^2}}}$ by multiplying a constant $h/\left( {{c}\sqrt \pi  } \right)$  and making change of the variable as $t \to t/{c}$. This results in
\begin{equation}\label{eq_2}
\frac{h}{{{c}\sqrt \pi  }}{e^{ - {{\left( {\frac{t}{{{c}}}} \right)}^2}}},
\end{equation}
where the ${c}$ is the fitting parameter. Figure 1 shows three variations of the sampling Gaussian function \eqref{eq_2} computed at $h = 0.25$,  ${c} = 0.15$ (black curve), $h = 0.25$,  ${c} = 0.2$ (blue curve) and $h = 0.25$,  ${c} = 0.25$ (red curve). We may attempt to approximate  $f\left( t \right)$ by replacing the sinc functions in equation \eqref{eq_1} with the Gaussian functions of kind \eqref{eq_2} as given by equation
$$
f\left( t \right) \approx \frac{h}{{{c}\sqrt \pi  }}\sum\limits_{n =  - N}^N {{e^{ - {{\left( {\frac{{t - {t_n}}}{{{c}}}} \right)}^2}}}f\left( {{t_n}} \right)}.
$$
Consequently, choosing the sampling points separated equidistantly with step $h$ between two adjacent points we obtain
\begin{equation}\label{eq_3}
f\left( t \right) \approx \frac{h}{{{c}\sqrt \pi  }}\sum\limits_{n =  - N}^N {{e^{ - {{\left( {\frac{{t - nh}}{{{c}}}} \right)}^2}}}f\left( {nh} \right)}.
\end{equation}

\begin{figure}[ht]
\begin{center}
\includegraphics[width=24pc]{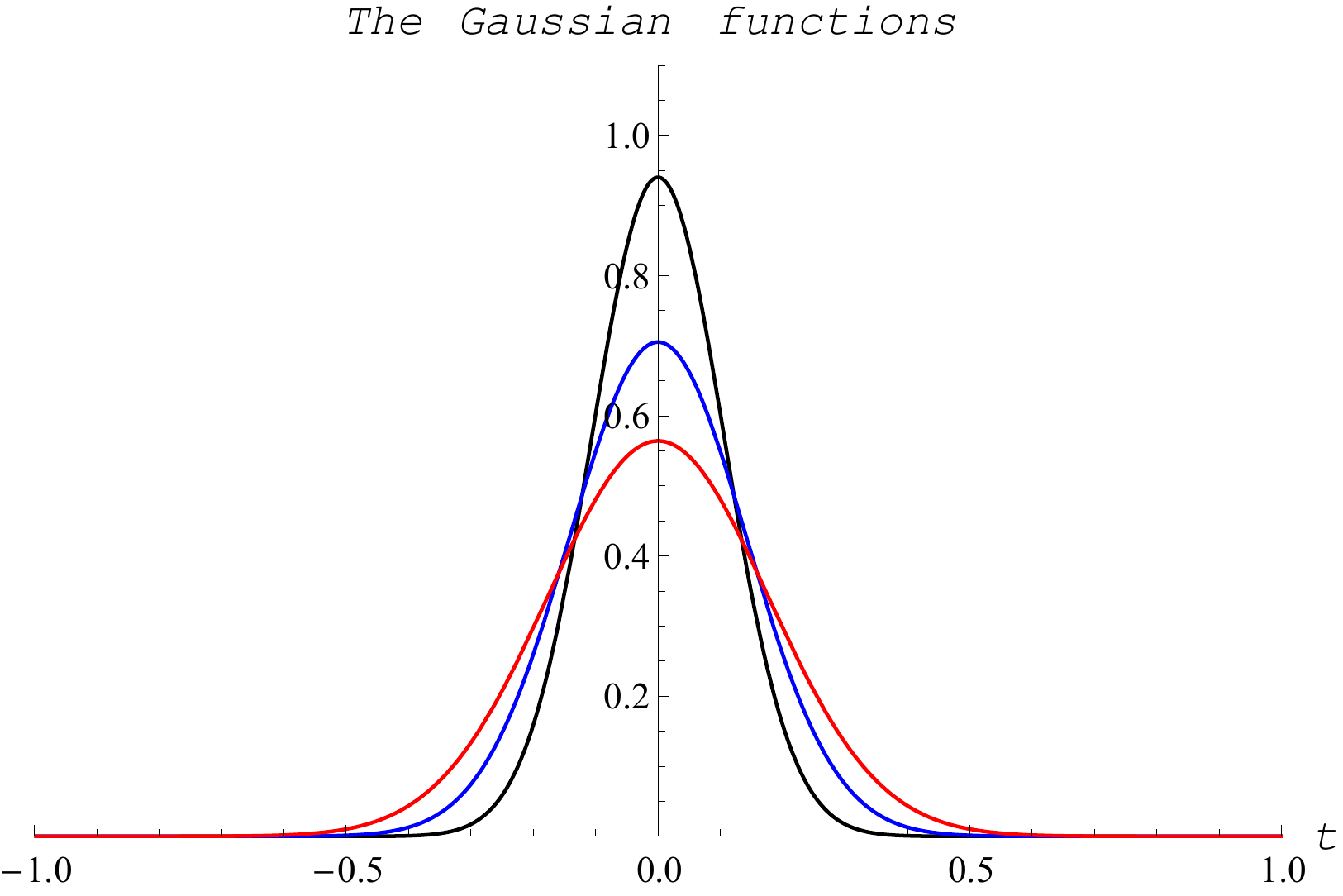}\hspace{2pc}%
\begin{minipage}[b]{28pc}
\vspace{0.3cm}
{\sffamily {\bf{Fig. 1.}} The Gaussian sampling functions computed at $h = 0.25$, ${c} = 0.15$ (black curve), $h = 0.25$, ${c} = 0.2$ (blue curve) and $h = 0.25$, ${c} = 0.25$ (red curve).}
\end{minipage}
\end{center}
\end{figure}

The application of the Gaussian function of kind \eqref{eq_2} can be justified by considering a simplest case when $f\left( t \right) = 1$. Figure 2 shows the corresponding approximations at $h = 0.25$, ${c} = 0.15$, $N = 10$ (blue curve) and $h = 0.25$, ${c} = 0.2$, $N = 10$ (red curve). The multiplication of the number of the sampling points to the step determines the range for approximation $\left( {2N + 1} \right)h = 5.25$. As a consequence, the approximated functions shown in Fig. 2 resemble the window function. In particular, these functions are approximately equal to $1$ within the range $t \in \left[ { - 5.25/2,5.25/2} \right]$ and to $0,$ otherwise. As we can see from this figure, at smaller ${c}$ the approximated function oscillates near the constant $1$ at the top of the curve stronger (blue curve). However, as ${c}$ increases the oscillation rapidly decreases (red curve) and practically vanishes at ${c} \mathbin{\lower.3ex\hbox{$\buildrel>\over
{\smash{\scriptstyle\sim}\vphantom{_x}}$}} h$. This signifies that if the change in curve between any two adjacent sampling points is relatively small, then any function can be approximated by sampling with the Gaussian function \eqref{eq_2}.

\begin{figure}[ht]
\begin{center}
\includegraphics[width=24pc]{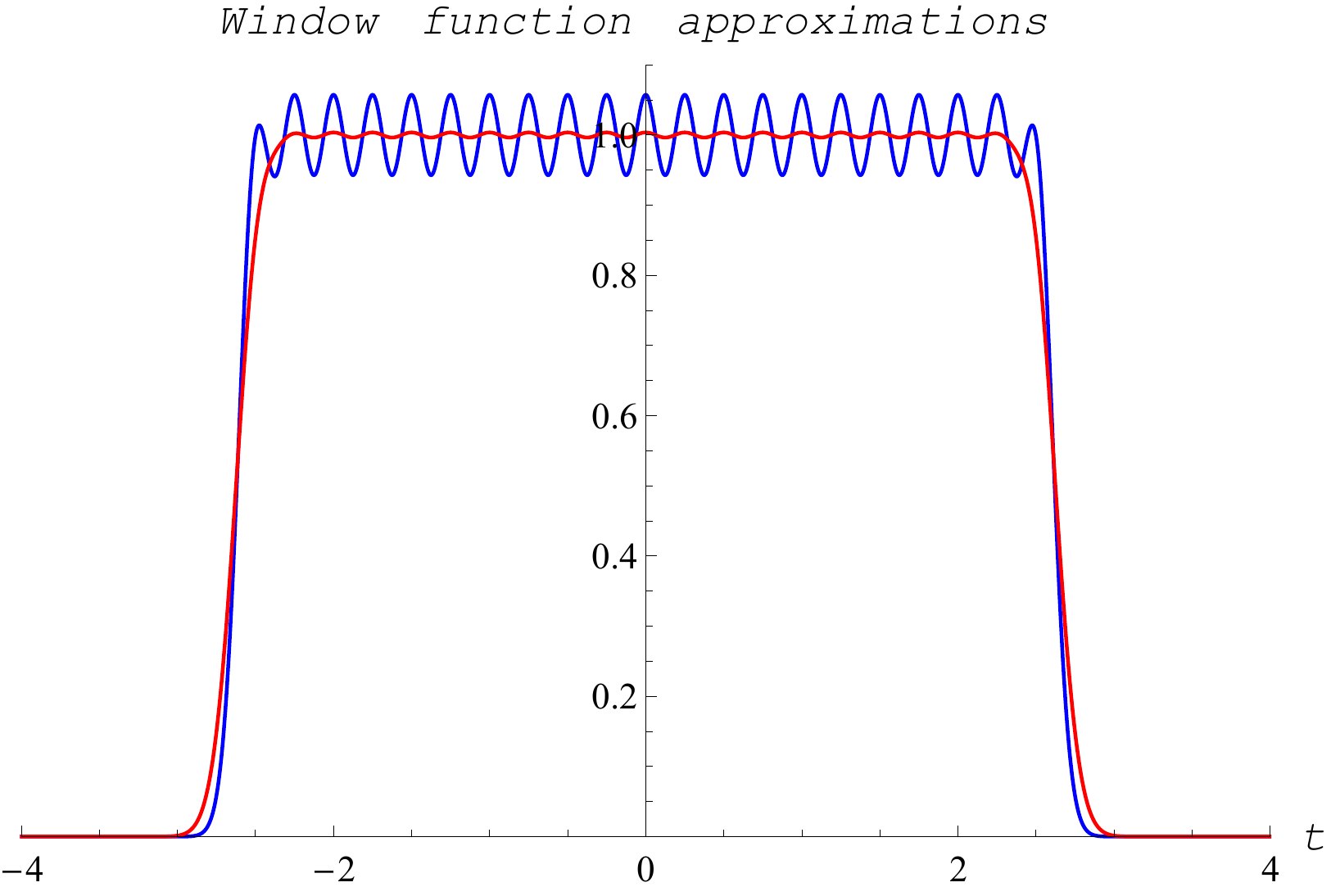}\hspace{2pc}%
\begin{minipage}[b]{28pc}
\vspace{0.3cm}
{\sffamily {\bf{Fig. 2.}}	The window function approximations computed at $h = 0.25$, ${c} = 0.15$, $N = 10$ (blue curve) and $h = 0.25$, ${c} = 0.2$, $N = 10$ (red curve).}
\end{minipage}
\end{center}
\end{figure}

In this work we introduce an application of the complex error function to the Fourier analysis. In particular, we show that the use of equation \eqref{eq_3} in the Fourier integration leads to a weighted sum of the complex error functions. Due to remarkable property of the complex error function this approach provides efficient computational methodology in the Fourier transform as a damping harmonic series.

\section{The complex error function}

The complex error function, also known as the Faddeeva function or the Kramp function, is defined as \cite{Faddeyeva1961, Gautschi1970, Abramowitz1972, Schreier1992}
\begin{equation}\label{eq_4}
w\left( z \right) = {e^{ - {z^2}}}\left( {1 + \frac{{2i}}{{\sqrt \pi  }}} \right)\int\limits_0^z {{e^{{u^2}}}du},
\end{equation}
where $z = x + iy$ is the complex argument. This function is a solution of the following differential equation \cite{Schreier1992}
$$
w'\left( z \right) + 2zw\left( z \right) = \frac{{2i}}{{\sqrt \pi  }}, \quad\quad w\left( 0 \right) = 1.
$$

The complex error function finds broad applications in many fields of Applied Mathematics \cite{Gautschi1970, Abramowitz1972, Weideman1994}, Physics \cite{Quine2002, Abrarov2015b, Quine2013, Berk2013, Borchert2003, Fried1961, Gordeyev1952} and Astronomy \cite{Emerson1996}. In Applied Mathematics it is closely related to the error function of complex argument \cite{Schreier1992}
\begin{equation}\label{eq_5}
w\left( {iz} \right) = {e^{ - {z^2}}}\left[ {1 - {\rm{erf}}\left( { - iz} \right)} \right] \Leftrightarrow {\rm{erf}}\left( z \right) = 1 - {e^{ - {z^2}}}w\left( {iz} \right),
\end{equation}
the normal distribution function \cite{Weisstein2003}
$$
\begin{aligned}
\Phi \left( z \right) &= \frac{1}{{\sqrt {2\pi } }}\int\limits_0^z {{e^{ - {u^2}/2}}du} \\
 &= \frac{1}{2}\left[ {1 - {e^{{ - z^2}/2}}w\left( {\frac{{iz}}{{\sqrt 2 }}} \right)} \right],
\end{aligned}
$$
the Fresnel integral \cite{Weisstein2003}
$$
\begin{aligned}
F\left( z \right) &= \int\limits_0^z {{e^{i\left( {\pi /2} \right){u^2}}}du} \\
 &= \left( {1 + i} \right)\left[ {1 - {e^{i\left( {\pi /2} \right){z^2}}}w\left( {\sqrt \pi  \left( {1 + i} \right)z/2} \right)} \right]/2
\end{aligned}
$$
and the Dawson\text{'}s integral \cite{Rybicki1989, Lether1998, Cody1970, McCabe1974, McKenna1984, Abrarov2015c}
$$
{\rm{daw}}\left( z \right) = {e^{ - {z^2}}}\int\limits_0^z {{e^{{u^2}}}du}  = \sqrt \pi  \frac{{ - {e^{ - {z^2}}} + w\left( z \right)}}{{2i}}.
$$

In Physics and Astronomy the complex error function is related to the Voigt function \cite{Schreier1992, Quine2002, Abrarov2015b, Quine2013, Berk2013, Borchert2003}
$$
K\left( {x,y} \right) = \left\{ 
\begin{aligned}
&\frac{y}{\pi }\int\limits_{ - \infty }^\infty  {\frac{{\exp \left( { - {u^2}} \right)}}{{{y^2} + {{\left( {x - u} \right)}^2}}}\,du}, \quad y \ne 0\\
&\exp \left( { - {x^2}} \right), \qquad\qquad\qquad y = 0
\end{aligned} \right.
$$
that describes the spectral behavior of the photon emitting or absorbing objects (photo-luminescent materials \cite{Borchert2003}, planetary atmosphere \cite{Schreier1992, Quine2002, Abrarov2015b, Quine2013, Berk2013}, celestial bodies \cite{Emerson1996} and so on). Specifically, the Voigt function represents the real part of the complex error function
$$
K\left( {x,y} \right) = {\mathop{\rm Re}\nolimits} \left[ {w\left( {x,y} \right)} \right], \quad\quad y \ge 0.
$$
Other functions that can be expressed in terms of the complex error function are the plasma dispersion function \cite{Fried1961}, the Gordeyev\text{'}s integral \cite{Gordeyev1952}, the rocket flight function \cite{Reichel1968} and the probability integral \cite{Weisstein2003}.

The complex error function \eqref{eq_4} can be represented alternatively as (see equation (3) in \cite{Srivastava1987} and \cite{Srivastava1992}, see also Appendix A in \cite{Abrarov2015c} for derivation)
$$
w\left( {x,y} \right) = \frac{1}{{\sqrt \pi  }}\int\limits_0^\infty  {\exp \left( { - {u^2}/4} \right)\exp \left( { - yu} \right)\exp \left( {ixu} \right)du}.
$$
Using the change of the variable as $u \to 2u$ in the integral above leads to \cite{Abrarov2014}
\begin{equation}\label{eq_6}
w\left( {x,y} \right) = \frac{2}{{\sqrt \pi  }}\int\limits_0^\infty  {{e^{ - {u^2}}}{e^{ - 2yu}}{e^{2ixu}}du}.
\end{equation}
Further, we will use this equation in derivation of the weighted sum.

 A rapid C/C++ implementation (RooFit package from CENR\text{'}s library) for computation of the complex error function with average accuracy $\sim{10^{ - 16}}$ has been reported in the recent work \cite{Karbach2014}.

\section{Derivations}

\subsection{Weighted sum}

There are several definitions for the Fourier transform. In this work we will use the following definitions \cite{Bracewell2000}
\begin{equation}\label{eq_7}
{\cal F}\left\{ {f\left( t \right)} \right\}\left( \nu  \right) = \int\limits_{ - \infty }^\infty  {f\left( t \right){e^{ - 2\pi i\nu t }}dt }\end{equation}
and
\begin{equation}\label{eq_8}
{{\cal F}^{ - 1}}\left\{ {F\left( \nu  \right)} \right\}\left( t \right) = \int\limits_{ - \infty }^\infty  {F\left( \nu  \right){e^{2\pi i\nu t }}d\nu}.
\end{equation}
Thus, the relationships between functions $f\left( t \right)$ and $F\left( \nu  \right)$ are performed by two operators ${\cal F}\left\{ {f\left( t \right)} \right\}\left( \nu  \right) = F\left( \nu  \right)$ and ${{\cal F}^{ - 1}}\left\{ {F\left( \nu  \right)} \right\}\left( t  \right) = f\left( t \right)$ corresponding to the forward and inverse Fourier transforms, respectively. In signal processing the arguments $t$ and $\nu $ in these reciprocally Fourier transformable functions $f\left( t \right)$ and $F\left( \nu  \right)$ are interpreted, accordingly, as time vs. frequency \cite{Hansen2014}.

We can find an approximation to the Fourier transform of the function $f\left( t \right)$ by substituting approximation \eqref{eq_3} into equation \eqref{eq_7}. This leads to
\begin{equation}\label{eq_9}
\begin{aligned}
F\left( \nu  \right) &\approx \frac{h}{{c\sqrt \pi  }}\int\limits_{ - \infty }^\infty  {\sum\limits_{n =  - N}^N {{e^{ - {{\left( {\frac{{t - nh}}{c}} \right)}^2}}}f\left( {nh} \right)} \,{e^{ - 2\pi i\nu t}}dt } \\
&= \frac{h}{{c\sqrt \pi  }}\sum\limits_{n =  - N}^N {\int\limits_{ - \infty }^\infty  {{e^{ - {{\left( {\frac{t}{c}} \right)}^2}}}{e^{2\frac{{nht}}{{{c^2}}}}}{e^{ - {{\left( {\frac{{nh}}{c}} \right)}^2}}}f\left( {nh} \right){e^{ - 2\pi i\nu t}}dt} }.
\end{aligned}
\end{equation}
Taking into account that
$$
\int\limits_{ - \infty }^\infty  {f\left( t \right){e^{ - 2\pi i\nu t}}dt}  = \int\limits_0^\infty  {f\left( t \right){e^{ - 2\pi i\nu t}}dt}  + \int\limits_0^\infty  {f\left( { - t} \right){e^{2\pi i\nu t}}dt},
$$
the equation \eqref{eq_9} can be rewritten as
\begin{equation}\label{eq_10}
\begin{aligned}
F\left( \nu  \right) \approx &\frac{h}{{c\sqrt \pi  }}\sum\limits_{n =  - N}^N {{e^{ - {{\left( {\frac{{nh}}{c}} \right)}^2}}}f\left( {nh} \right)\int\limits_0^\infty  {{e^{ - {{\left( {\frac{t}{c}} \right)}^2}}}{e^{2\frac{{nht}}{{{c^2}}}}}{e^{ - 2\pi i\nu t}}dt} } \\
& + \frac{h}{{c\sqrt \pi  }}\sum\limits_{n =  - N}^N {{e^{ - {{\left( {\frac{{nh}}{c}} \right)}^2}}}f\left( { - nh} \right)\int\limits_0^\infty  {{e^{ - {{\left( {\frac{t}{c}} \right)}^2}}}{e^{2\frac{{nht}}{{{c^2}}}}}{e^{2\pi i\nu t}}dt} }.
\end{aligned}
\end{equation}
Since change of the variable $t/c \to t$ leads to
\[
\begin{aligned}
\int\limits_0^\infty  \exp \left( { - \frac{{{t^2}}}{{{c^2}}}} \right)\exp \left( {2\frac{{nht}}{{{c^2}}}} \right)&\exp \left( { \pm 2\pi it\nu } \right)dt  \\
&= c\int\limits_0^\infty  {\exp \left( { - {t^2}} \right)\exp \left( {2\frac{{nht}}{c}} \right)\exp \left( { \pm 2\pi ic\nu t} \right)dt},
\end{aligned}
\]
the equation \eqref{eq_10} can be rearranged as
\[
\begin{aligned}
F\left( \nu  \right) \approx &\frac{h}{{\sqrt \pi  }}\sum\limits_{n =  - N}^N {{e^{ - {{\left( {\frac{{nh}}{c}} \right)}^2}}}f\left( {nh} \right)\int\limits_0^\infty  {{e^{ - {t^2}}}{e^{2\frac{{nh}}{c}t}}{e^{ - 2\pi ci\nu t}}dt} }  \\
&+ \frac{h}{{\sqrt \pi  }}\sum\limits_{n =  - N}^N {{e^{ - {{\left( {\frac{{nh}}{c}} \right)}^2}}}f\left( { - nh} \right)\int\limits_0^\infty  {{e^{ - {t^2}}}{e^{2\frac{{nh}}{c}t}}{e^{2\pi ci\nu t}}dt} }.
\end{aligned}
\]
Comparing this equation with equation \eqref{eq_6} yields
\begin{equation}\label{eq_11}
\begin{aligned}
F\left( \nu  \right) \approx &\frac{h}{2}\sum\limits_{n =  - N}^N {{e^{ - {{\left( {\frac{{nh}}{c}} \right)}^2}}}f\left( {nh} \right)w\left( { - \pi c\nu , - nh/c} \right)}  \\
&+ \frac{h}{2}\sum\limits_{n =  - N}^N {{e^{ - {{\left( {\frac{{nh}}{c}} \right)}^2}}}f\left( { - nh} \right)w\left( {\pi c\nu , - nh/c} \right)}.
\end{aligned}
\end{equation}
Lastly, defining the constants ${\alpha _n} = h {e^{ - {{\left( {\frac{{nh}}{c}} \right)}^2}}}f\left( {nh} \right) /2$ we can rewrite approximation \eqref{eq_11} in a more compact form as a weighted sum
\begin{equation}\label{eq_12}
\begin{aligned}
F\left( \nu  \right) & = {\cal F} \left\{ f \left( t \right) \right\} \left( \nu \right) \\
&\approx \sum\limits_{n =  - N}^{n = N} {\left[ {{\alpha _n}w\left( { - \pi c\nu , - nh/c} \right) + {\alpha _{ - n}}w\left( {\pi c\nu , - nh/c} \right)} \right]}.
\end{aligned}
\end{equation}

The equation \eqref{eq_12} can be expressed through any functions that have been considered in the section $2$ above. For example, using the identity \eqref{eq_5} after trivial rearrangements we get the Fourier transform in terms of error functions of complex argument
\[
\begin{aligned}
F\left( \nu  \right) \approx \sum\limits_{n =  - N}^{n = N} &\left\{ {\alpha _n}{e^{ - {{\left( {\pi c\nu  + inh/c} \right)}^2}}}\left[ {1 + {\rm{erf}}\left( {nh/c - i\pi c\nu } \right)} \right] \right. \\
&\left. \qquad \, + {\alpha _{ - n}}{e^{ - {{\left( {\pi c\nu  - inh/c} \right)}^2}}}\left[ {1 + {\rm{erf}}\left( {nh/c + i\pi c\nu} \right)} \right] \right\} .
\end{aligned}
\]
It can also be shown that the approach based on a weighted sum of the complex error functions can be generalized to the Laplace transform.

The derivation of approximation for the inverse Fourier transform is straightforward now. Comparing integrals \eqref{eq_7}, \eqref{eq_8} and using approximation \eqref{eq_12} we immediately obtain
\begin{equation}\label{eq_13}
\begin{aligned}
f\left( t \right) &= {{\cal F}^{ - 1}}\left\{ {F\left( \nu  \right)} \right\}\left( t \right) \\
&\approx \sum\limits_{n =  - N}^{n = N} {\left[ {{\alpha _n^*}w\left( {\pi ct, - nh/c} \right) + {\alpha _{ - n}^*}w\left( { - \pi ct, - nh/c} \right)} \right]},
\end{aligned}
\end{equation}
where the coefficients are calculated as ${\alpha _n^*} = h{e^{ - {{\left( {\frac{{nh}}{c}} \right)}^2}}}F\left( {nh} \right) /2$.

Since the Gaussian function \eqref{eq_2} rapidly decreases with increasing $\left| t \right|$, only few terms with negative index $n$ actually contribute to shape the curve $f\left( t \right)$ along the positive $t$-axis where the Fourier integration takes place according to equation \eqref{eq_10}. As a result, it is sufficient to take into consideration only, say, first three terms $h{e^{ - {{\left[ {\left( {t + 3h} \right)/c} \right]}^2}}}/\left( {c\sqrt \pi  } \right)$, $h{e^{ - {{\left[ {\left( {t + 2h} \right)/c} \right]}^2}}}/\left( {c\sqrt \pi  } \right)$ and $h{e^{ - {{\left[ {\left( {t + h} \right)/c} \right]}^2}}}/\left( {c\sqrt \pi  } \right)$ with negative index $n$. Thus, the number of the sampling points is reduced by almost factor of $2$ from $2N + 1$ to $N + 4$. Consequently, the equations \eqref{eq_12} and \eqref{eq_13} can be simplified as
\begin{equation}\label{eq_14}
\begin{aligned}
F\left( \nu  \right) &= {\cal F}\left\{ {f\left( t \right)} \right\}\left( \nu  \right) \\
&\approx \sum\limits_{n =  - 3}^{n = N} {\left[ {{\alpha _n}w\left( { - \pi c\nu , - nh/c} \right) + {\alpha _{ - n}}w\left( {\pi c\nu , - nh/c} \right)} \right]}
\end{aligned}
\end{equation}
and
\begin{equation}\label{eq_15}
\begin{aligned}
f\left( t \right) &= {{\cal F}^{ - 1}}\left\{ {F\left( \nu  \right)} \right\}\left( t \right) \\
&\approx \sum\limits_{n =  - 3}^{n = N} {\left[ {{\alpha _n^*}w\left( {\pi ct, - nh/c} \right) + {\alpha _{ - n}^*}w\left( { - \pi ct, - nh/c} \right)} \right]}.
\end{aligned}
\end{equation}

It should be noted that a {\it{real-time}} computation of the complex error functions in equations \eqref{eq_14} and \eqref{eq_15} may not be fast enough for the wavelets with extended effective length. However, this problem can be easily resolved by using the following precomputed values
$$
\begin{aligned}
&\beta_{n,1} = {\frac{h}{2}}{{e^{ - {{\left( {\frac{{nh}}{c}} \right)}^2}}}w\left( { - \pi c\nu , - nh/c} \right)},\\
&\beta_{n,2} ={\frac{h}{2}}{{e^{ - {{\left( {\frac{{nh}}{c}} \right)}^2}}}w\left( {\pi c\nu , - nh/c} \right)}
\end{aligned}
$$
and
$$
\begin{aligned}
&\beta_{n,1}^* = {\frac{h}{2}}{{e^{ - {{\left( {\frac{{nh}}{c}} \right)}^2}}}w\left( {\pi ct , - nh/c} \right)},\\
&\beta_{n,2}^* ={\frac{h}{2}}{{e^{ - {{\left( {\frac{{nh}}{c}} \right)}^2}}}w\left( {- \pi ct , - nh/c} \right)}
\end{aligned}
$$
that can be stored as datasets in a computer memory. In particular, the application of these {\it{ready-made}} precomputed values enables us to simplify considerably the approximations \eqref{eq_14} and \eqref{eq_15} for the forward and inverse Fourier transforms as given by
\begin{equation}\label{eq_16}
F\left( \nu  \right) = {\cal F}\left\{ {f\left( t \right)} \right\}\left( \nu  \right) \approx \sum\limits_{n =  - 3}^N {f\left( {nh} \right) \beta_{n,1}} + \sum\limits_{n =  - 3}^N {f\left( { - nh} \right) \beta_{n,2}}
\end{equation}
and
\begin{equation}\label{eq_17}
f\left( t  \right) = {{\cal F}^{ - 1}}\left\{ {F\left( \nu  \right)} \right\}\left( t \right) \approx \sum\limits_{n =  - 3}^N {F\left( {nh} \right) \beta_{n,1}^*} + \sum\limits_{n =  - 3}^N {F\left( { - nh} \right) \beta_{n,2}^*} \, ,
\end{equation}
accordingly.

\subsection{Damping harmonic series}

As an application of the precomputed values $\beta_{n,1}$, $\beta_{n,2}$ and $\beta_{n,1}^*$, $\beta_{n,2}^*$ in approximations \eqref{eq_16} and \eqref{eq_17} helps us to avoid a {\it{real-time}} computation of the complex error functions, this technique significantly accelerates an algorithmic implementation. However, there is an alternative way for rapid performance in the numerical forward and inverse Fourier transforms. Specifically, using a remarkable property of the complex error function \cite{McKenna1984, Abrarov2014}
\begin{equation}\label{eq_18}
w\left( z \right) = 2{e^{ - {z^2}}} - w\left( { - z} \right)
\end{equation}
the approximations \eqref{eq_12} and \eqref{eq_13} can be significantly simplified and expressed in terms of a damping harmonic series. 

Let us rewrite equation \eqref{eq_12} as
\[
\begin{aligned}
F\left( \nu  \right) \approx &\left[ {{\alpha _0}w\left( { - \pi c\nu ,0} \right) + {\alpha _0}w\left( {\pi c\nu ,0} \right)} \right]\\
 &+ \sum\limits_{n = 1}^N {\left[ {{\alpha _n}w\left( { - \pi c\nu , - nh/c} \right) + {\alpha _{ - n}}w\left( {\pi c\nu , - nh/c} \right)} \right]} \\
 &+ \sum\limits_{n =  - N}^{ - 1} {\left[ {{\alpha _n}w\left( { - \pi c\nu , - nh/c} \right) + {\alpha _{ - n}}w\left( {\pi c\nu , - nh/c} \right)} \right]}
\end{aligned}
\]
or
\begin{equation}
\begin{aligned}\label{eq_19}
F\left( \nu  \right) \approx &\left[ {{\alpha _0}w\left( { - \pi c\nu ,0} \right) + {\alpha _0}w\left( {\pi c\nu ,0} \right)} \right]\\
 &+ \sum\limits_{n = 1}^N {\left[ {{\alpha _n}w\left( { - \pi c\nu , - nh/c} \right) + {\alpha _{ - n}}w\left( {\pi c\nu , - nh/c} \right)} \right]} \\
 &+ \sum\limits_{n = 1}^N {\left[ {{\alpha _{ - n}}w\left( { - \pi c\nu ,nh/c} \right) + {\alpha _n}w\left( {\pi c\nu ,nh/c} \right)} \right]} .
\end{aligned}
\end{equation}

The identity \eqref{eq_18} can be expressed through arguments $x$ and $y$ as
$$
w\left( {x,y} \right) = 2{e^{ - {{\left( {x + iy} \right)}^2}}} - w\left( { - x, - y} \right) = 2{e^{{y^2} - {x^2}}}{e^{ - 2ixy}} - w\left( { - x, - y} \right).
$$
This leads to
$$
w\left( { \pm \pi c\nu , - nh/c} \right) = 2{e^{{{\left( {nh/c} \right)}^2} - {{\left( {\pi c\nu } \right)}^2}}}{e^{ \pm 2\pi i\nu nh}} - w\left( { \mp \pi c\nu ,nh/c} \right).
$$
Therefore, we have the following identity
\[
\begin{aligned}
\frac{h}{2}&\sum\limits_{n = 1}^N {\left[ {{\alpha _n}w\left( { - \pi c\nu , - nh/c} \right) + {\alpha _{ - n}}w\left( {\pi c\nu , - nh/c} \right)} \right]} \\
 &  = h\sum\limits_{n = 1}^N {\left[ {f\left( {nh} \right){e^{ - {{\left( {\pi c\nu } \right)}^2} - 2\pi i\nu nh}} + f\left( { - nh} \right){e^{ - {{\left( {\pi c\nu } \right)}^2} + 2\pi i\nu nh}}} \right]} \\
 & \,\,\,\,\, - \frac{h}{2}\sum\limits_{n = 1}^N {\left[ {{\alpha _n}w\left( {\pi c\nu ,nh/c} \right) + {\alpha _{ - n}}w\left( { - \pi c\nu ,nh/c} \right)} \right],} 
\end{aligned}
\]
where $hf\left( {nh} \right) = 2{\alpha _n}{e^{{{\left( {\frac{{nh}}{c}} \right)}^2}}}$.
Substituting this identity into approximation \eqref{eq_19} yields
$$
\begin{aligned}
F\left( \nu  \right) \approx &\frac{h}{2}f\left( 0 \right)\left[ {w\left( { - \pi c\nu ,0} \right) + w\left( {\pi c\nu ,0} \right)} \right] \\
&+ h{e^{ - {{\left( {\pi c\nu } \right)}^2}}}\sum\limits_{n = 1}^N {\left[ {f\left( {nh} \right){e^{ - 2\pi i\nu nh}} + f\left( { - nh} \right){e^{2\pi i\nu nh}}} \right]}.
\end{aligned}
$$
Since (see \cite{Abrarov2015b} for details)
$$
w\left( { - x,0} \right) + w\left( {x,0} \right) = 2{e^{ - {x^2}}},
$$
this approximation can be expressed as
$$
\begin{aligned}
F\left( \nu  \right) \approx &hf\left( 0 \right){e^{ - {{\left( {\pi c\nu } \right)}^2}}} \\
&+ h{e^{ - {{\left( {\pi c\nu } \right)}^2}}}\sum\limits_{n = 1}^N {\left[ {f\left( {nh} \right){e^{ - 2\pi i\nu nh}} + f\left( { - nh} \right){e^{2\pi i\nu nh}}} \right]}
\end{aligned}
$$
or
\begin{equation}\label{eq_20}
F\left( \nu  \right) = {\cal F}\left\{ {f\left( t \right)} \right\}\left( \nu  \right) \approx h{e^{ - {{\left( {\pi c\nu } \right)}^2}}}\sum\limits_{n =  - N}^N {f\left( {nh} \right){e^{ - 2\pi i\nu nh}}}.
\end{equation}

The inverse Fourier transform can be found from approximation \eqref{eq_13} by using the same derivation procedure
\begin{equation}\label{eq_21}
f\left( t \right) = {{\cal F}^{ - 1}}\left\{ {F\left( \nu  \right)} \right\}\left( t \right) \approx h{e^{ - {{\left( {\pi ct} \right)}^2}}}\sum\limits_{n =  - N}^N {F\left( {nh} \right){e^{2\pi itnh}}}.
\end{equation}

Although the approximations \eqref{eq_20} and \eqref{eq_21} consist of harmonic functions ${e^{ - 2\pi i\nu nh}} = cos\left({2\pi \nu nh }\right) - i \, sin\left({2\pi \nu nh }\right)$ and ${e^{ 2\pi itnh}} = cos\left({2\pi tnh }\right) + i \, sin\left({2\pi tnh }\right)$, the resultant forward and inverse Fourier transformed wavelets remain, nevertheless, non-periodic due to presence of the damping exponential functions ${e^{ - {{\left( {\pi c\nu} \right)}^2}}}$ and ${e^{ - {{\left( {\pi ct} \right)}^2}}}$. Consequently, unlike the conventional discrete Fourier transform the approximations \eqref{eq_20} and \eqref{eq_21} have no distortion that may appear as a result of aliasing.

\section{Sample of computation}

If a function $f \left( t \right)$ can be separated into the even and odd parts as
$$
f\left( t \right) = \underbrace {f^+\left( t \right)}_{even} + \underbrace {f^-\left( t \right)}_{odd},
$$
where $f^+ \left( - t \right) = f^+ \left( t \right)$ and $f^- \left( - t \right) = - f^- \left( t \right)$, then according to definition \eqref{eq_7} the forward Fourier transforms are given by
\begin{equation}\label{eq_22}
{\cal F}\left\{ {f^+\left( t \right)} \right\}\left( \nu  \right) = \int\limits_{ - \infty }^\infty  {f^+\left( t \right)\cos \left( {2\pi \nu t} \right)dt }
\end{equation}
and
\begin{equation}\label{eq_23}
{\cal F}\left\{ {f^-\left( t \right)} \right\}\left( \nu  \right) =  - i\int\limits_{ - \infty }^\infty  {f^-\left( t \right)\sin \left( {2\pi \nu t} \right)dt },
\end{equation}
respectively. Therefore, using approximations \eqref{eq_20} and \eqref{eq_21} for the forward Fourier transforms over the even $f^+\left( t \right)$ and odd $f^-\left( t \right)$ parts we obtain
\begin{equation}\label{eq_24}
G\left( \nu  \right) \approx h{e^{ - {{\left( {\pi c\nu } \right)}^2}}}\sum\limits_{n =  - N}^N {f^+\left( {nh} \right)\cos \left( {2\pi \nu nh} \right)}
\end{equation}
and
\begin{equation}\label{eq_25}
H\left( \nu  \right) \approx  - ih{e^{ - {{\left( {\pi c\nu } \right)}^2}}}\sum\limits_{n =  - N}^N {f^-\left( {nh} \right)\sin \left( {2\pi \nu nh} \right)},
\end{equation}
respectively. Since
\[
\begin{aligned}
\sum\limits_{n =  -N}^{-1} {f^+\left( {nh} \right)\cos \left( {2\pi \nu nh} \right)} &= \sum\limits_{n =  1}^N {f^+\left( {nh} \right)\cos \left( {2\pi \nu nh} \right)},\\
{f^+\left( {0h} \right)\cos \left( {2\pi \nu 0h} \right)} &= {f^+\left( {0} \right)}
\end{aligned}
\]
and
\[
\begin{aligned}
\sum\limits_{n =  -N}^{-1} {f^-\left( {nh} \right)\sin \left( {2\pi \nu nh} \right)} &= \sum\limits_{n =  1}^N {f^-\left( {nh} \right)\sin \left( {2\pi \nu nh} \right)},\\
{f^-\left( {0h} \right)\sin \left( {2\pi \nu 0 h} \right)} &= 0
\end{aligned}
\]
from equations \eqref{eq_24} and \eqref{eq_25} it follows that
\begin{equation}\label{eq_26}
\begin{aligned}
G\left( \nu  \right) &= {\cal F}\left\{ {f^+\left( t \right)} \right\}\left( \nu  \right) \\
&\approx 2h{e^{ - {{\left( {\pi c\nu } \right)}^2}}} \left(f^+\left( 0 \right)/2 + \sum\limits_{n =  1}^N {f^+\left( {nh} \right)\cos \left( {2\pi \nu nh} \right)} \right)
\end{aligned}
\end{equation}
and
\begin{equation}\label{eq_27}
\begin{aligned}
H\left( \nu  \right) &= {\cal F}\left\{ {f^-\left( t \right)} \right\}\left( \nu  \right) \\
&\approx  - 2ih{e^{ - {{\left( {\pi c\nu } \right)}^2}}}\sum\limits_{n =  1}^N {f^-\left( {nh} \right)\sin \left( {2\pi \nu nh} \right)}.
\end{aligned}
\end{equation}

\begin{figure}[ht]
\begin{center}
\includegraphics[width=24pc]{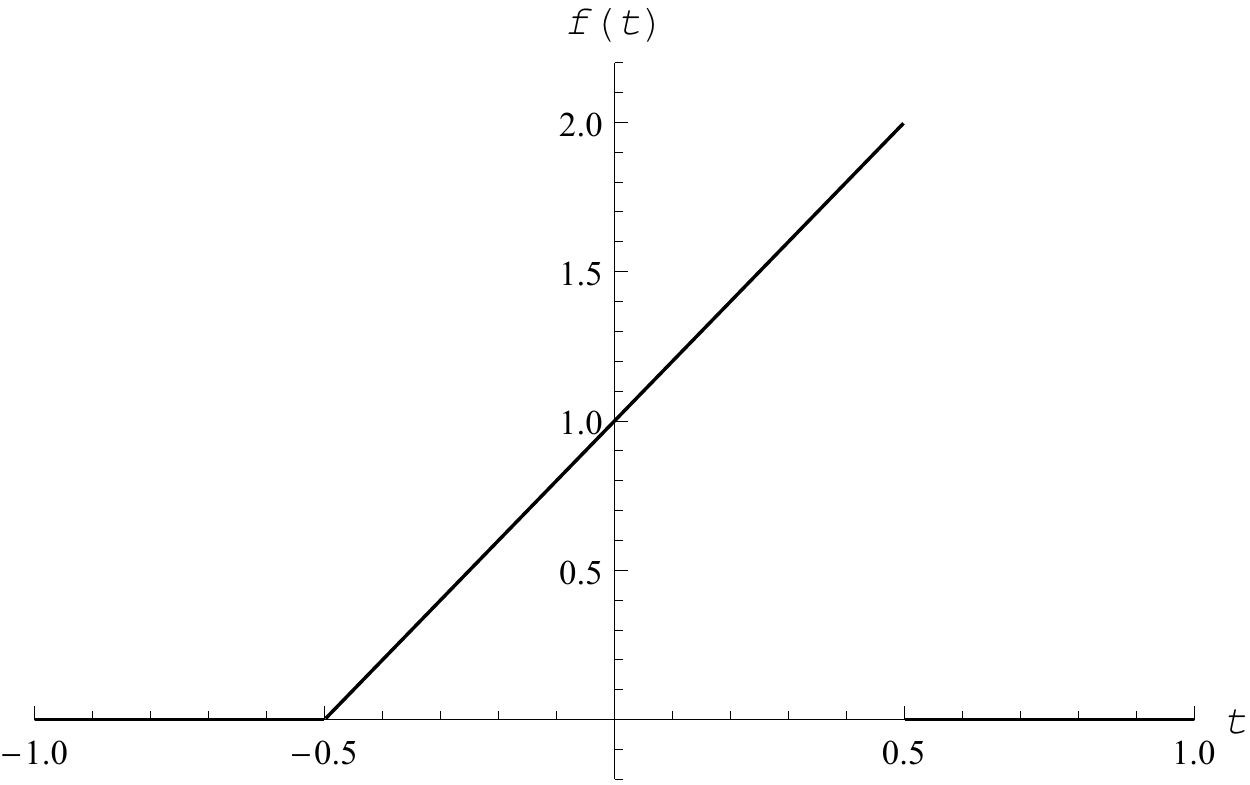}\hspace{2pc}%
\begin{minipage}[b]{28pc}
\vspace{0.3cm}
{\sffamily {\bf{Fig. 3.}} The function $f\left( t \right)$ defined according to equation \eqref{eq_28}.}
\end{minipage}
\end{center}
\end{figure}

Consider the following example
\begin{equation}\label{eq_28}
f\left( t \right) = \left\{ \begin{aligned}
2t+1, &  \qquad - 1/2 \le t \le 1/2\\
0, & \qquad {\rm{otherwise}}.
\end{aligned} 
\right.
\end{equation}
The equation \eqref{eq_28} can be separated into the even
\begin{equation}\label{eq_29}
f^+\left( t \right) = \left\{ \begin{aligned}
1, &  \qquad - 1/2 \le t \le 1/2\\
0, & \qquad {\rm{otherwise}}
\end{aligned} 
\right.
\end{equation}
and odd 
\begin{equation}\label{eq_30}
f^-\left( t \right) = \left\{ \begin{aligned}
2t, &  \qquad - 1/2 \le t \le 1/2\\
0, & \qquad {\rm{otherwise}}
\end{aligned}
\right.
\end{equation}
parts. Figure 3 shows that the function $f \left( t \right)$ while Fig. 4 depicts its even $f^+ \left( t \right)$ and odd $f^- \left( t \right)$ parts by blue and red curves, respectively.

From Fig. 4 we can see that the functions $f^+\left( t \right)$ and $f^-\left( t \right)$ are the wavelets (or the pulses) localized within the region $t \in \left[ { - 1/2,1/2} \right]$. Beyond this region the functions $f^+\left( t \right)$ and $f^-\left( t \right)$ are equal to zero. Therefore, the effective length of these wavelets is $\Delta t = 2 \times 1/2 = 1$. Since according to approximation \eqref{eq_3} we applied $2N + 1$ sampling points, the step $h$ between two adjacent sampling points can be determined from the formula $\left( {2N + 1} \right)h = \Delta t$. Thus, by choosing $N = 50$ and $N = 300$ we can find the corresponding steps to be $h = 0.0099$ and $h = 0.00166389$, respectively.

\begin{figure}[ht]
\begin{center}
\includegraphics[width=24pc]{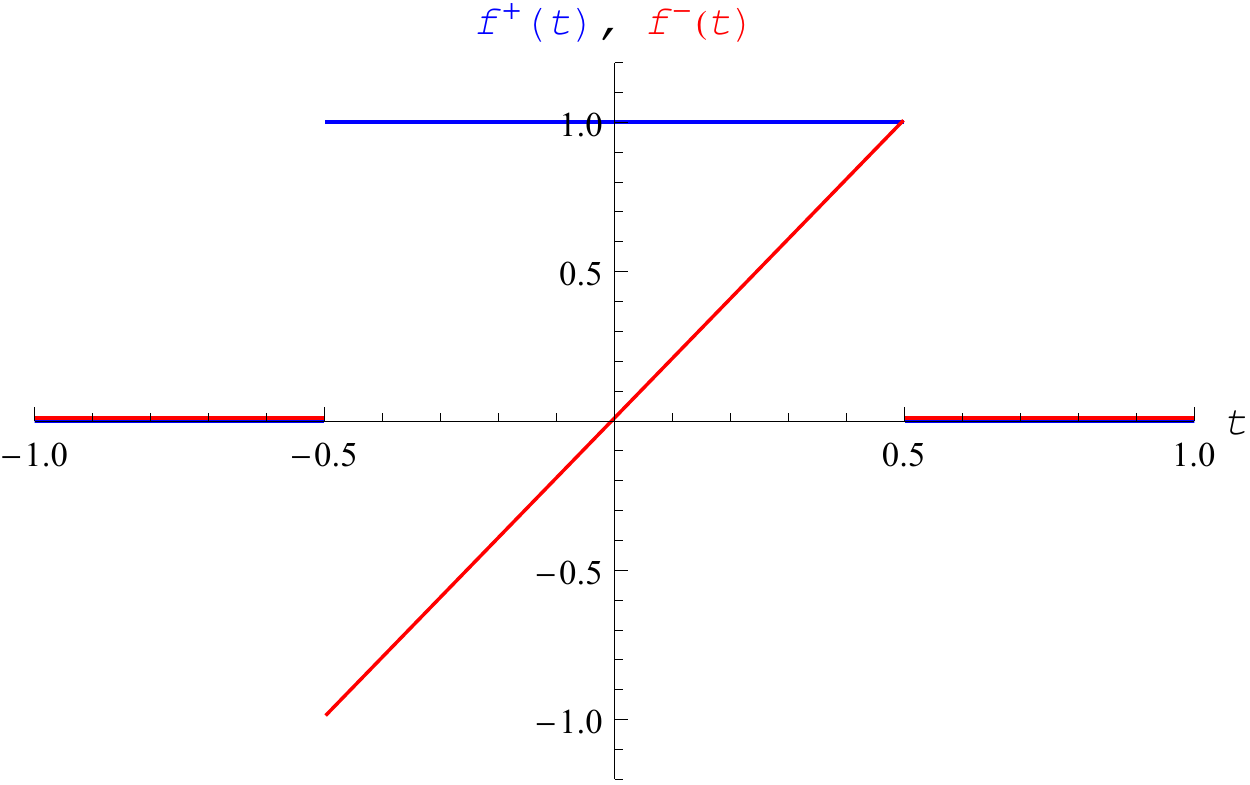}\hspace{2pc}%
\begin{minipage}[b]{28pc}
\vspace{0.3cm}
{\sffamily {\bf{Fig. 4.}} The even $f^+\left( t \right)$ (blue) and odd $f^-\left( t \right)$ (red) parts of the function $f \left( t \right)$.}
\end{minipage}
\end{center}
\end{figure}

As the function $f^- \left( t \right)$ is odd, its Fourier transform $H\left( \nu  \right)$ is purely imaginary according to equation \eqref{eq_23}. Figure 5 depicts $G\left( \nu  \right)$ and $\mathop{\rm Im} \left[ H\left( \nu  \right) \right]$ corresponding to the functions $f^+ \left( t \right)$ and $f^- \left( t \right)$ by blue and red curves, respectively. The curves for the functions $G\left( \nu  \right)$ and $\mathop{\rm Im} \left[ H\left( \nu  \right) \right]$ are obtained numerically by using equations \eqref{eq_26} and \eqref{eq_27} at $N = 50$, $h = c = 0.0099$.

\begin{figure}[ht]
\begin{center}
\includegraphics[width=24pc]{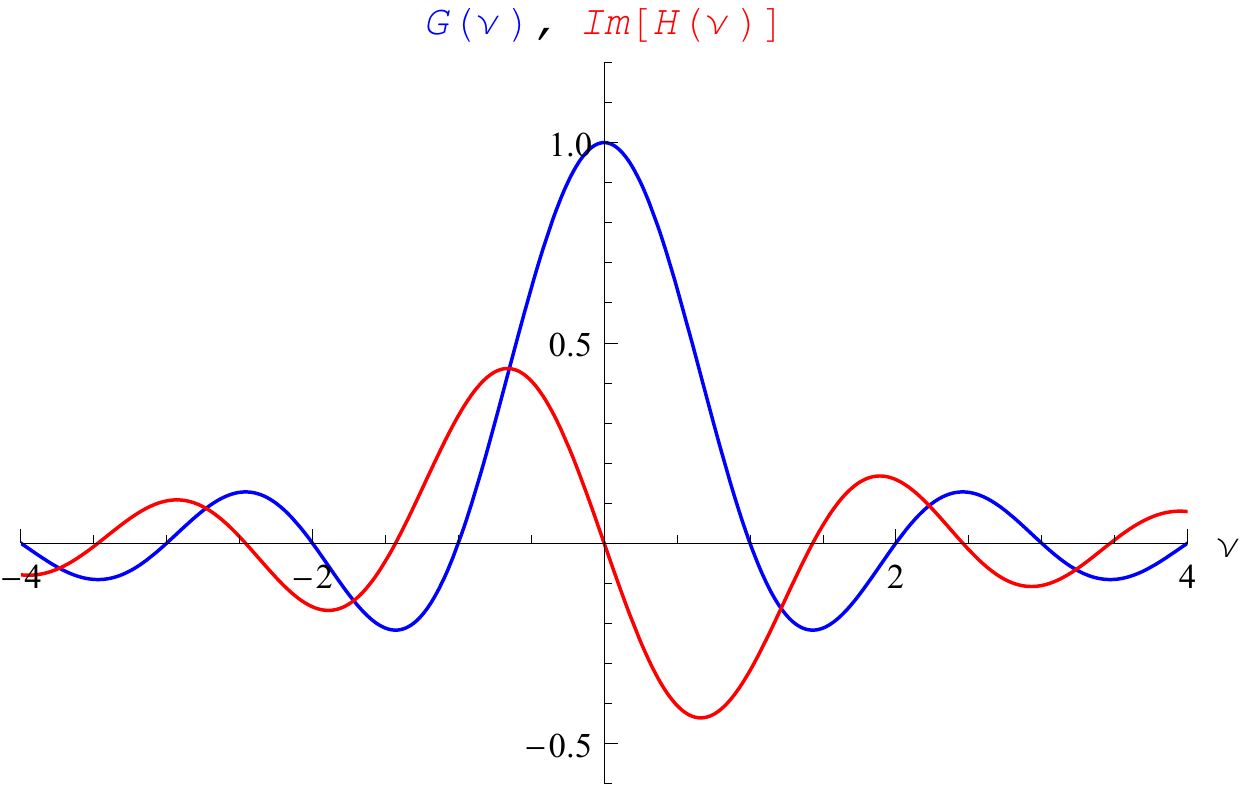}\hspace{2pc}%
\begin{minipage}[b]{28pc}
\vspace{0.3cm}
{\sffamily {\bf{Fig. 5.}} The Fourier transform functions corresponding to the even (blue) and odd (red) parts of the equation \eqref{eq_28}.}
\end{minipage}
\end{center}
\end{figure}

The Fourier transforms for the even $f^+ \left( t \right)$ and odd $f^- \left( t \right)$ parts of the function $f \left( t \right)$ can be found analytically. In particular, substituting the equations \eqref{eq_29}, \eqref{eq_30} into approximations \eqref{eq_22}, \eqref{eq_23} and considering the fact that these wavelets are not zero-valued only at $t \in \left[ { - 1/2,1/2} \right]$ we get
$$
G\left( \nu  \right) = \int\limits_{ - 1/2}^{1/2} {{\cos \left( {2\pi \nu t} \right)}dt = } \,\,{\rm{sinc}}\left( {\pi \nu } \right)
$$
and
$$
H\left( \nu  \right) = -i\int\limits_{ - 1/2}^{1/2} {2t \, {\sin\left( {2\pi \nu t} \right)}dt  = } \,\,i\frac{{\pi \nu \cos \left( {\pi \nu } \right) - \sin \left( {\pi \nu } \right)}}{{{\pi ^2}{\nu ^2}}}.
$$
Therefore, it is convenient to define the differences by using these functions
$$
\begin{aligned}
{\Delta _{{\mathop{\rm Re}\nolimits} }} &= G\left( \nu  \right) - 2h{e^{ - {{\left( {\pi c\nu } \right)}^2}}} \left(f^+\left( 0 \right)/2 + \sum\limits_{n =  1}^N {f^+\left( {nh} \right)\cos \left( {2\pi \nu nh} \right)} \right) \\
 &= {\rm{sinc}}\left( {\pi \nu } \right) -  2h{e^{ - {{\left( {\pi c\nu } \right)}^2}}} \left(f^+\left( 0 \right)/2 + \sum\limits_{n =  1}^N {f^+\left( {nh} \right)\cos \left( {2\pi \nu nh} \right)} \right)
\end{aligned}
$$
and
$$
\begin{aligned}
{\Delta _{{\mathop{\rm Im}\nolimits} }} &= {\mathop{\rm Im}\nolimits} \left[ {H\left( \nu  \right)} \right] - \left\{
- 2h{e^{ - {{\left( {\pi c\nu } \right)}^2}}}\sum\limits_{n =  1}^N {f^-\left( {nh} \right)\sin \left( {2\pi \nu nh} \right)}
\right\}\\
 &= \frac{{\pi \nu \cos \left( {\pi \nu } \right) - \sin \left( {\pi \nu } \right)}}{{{\pi ^2}{\nu ^2}}} - \left\{
- 2h{e^{ - {{\left( {\pi c\nu } \right)}^2}}}\sum\limits_{n =  1}^N {f^-\left( {nh} \right)\sin \left( {2\pi \nu nh} \right)} \right\}.
\end{aligned}
$$

\begin{figure}[ht]
\begin{center}
\includegraphics[width=24pc]{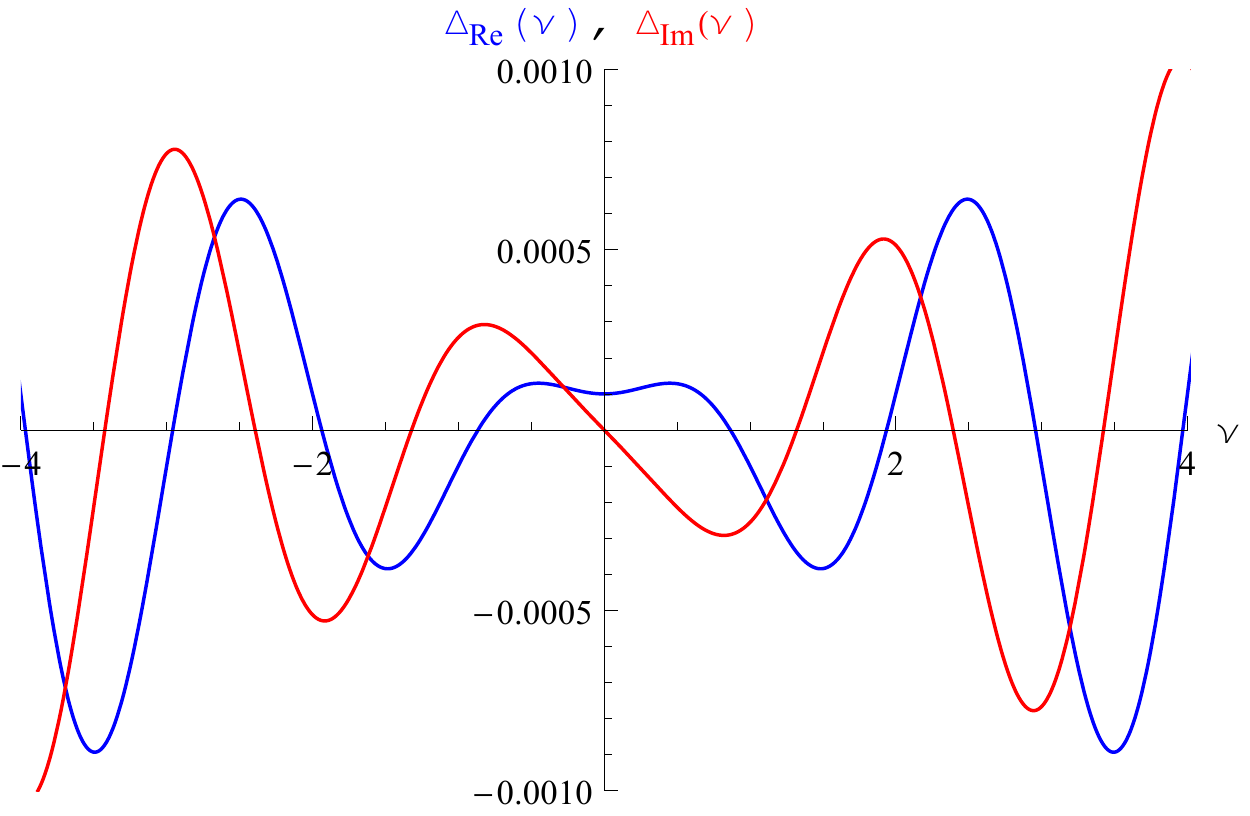}\hspace{2pc}%
\begin{minipage}[b]{28pc}
\vspace{0.3cm}
{\sffamily {\bf{Fig. 6.}} The differences ${\Delta _{{\mathop{\rm Re}\nolimits} }}$ (blue curve) and ${\Delta _{{\mathop{\rm Im}\nolimits} }}$ (red curve) computed at $N = 50$, $h = c = 0.0099$, respectively.}
\end{minipage}
\end{center}
\end{figure}

Figure 6 illustrates the differences ${\Delta _{{\mathop{\rm Re}\nolimits} }}$ (blue curve) and ${\Delta _{{\mathop{\rm Im}\nolimits} }}$ (red curve) computed at $N = 50$, $h = c = 0.0099$. As we can see from Fig. 6, the differences ${\Delta _{{\mathop{\rm Re}\nolimits} }}$ and ${\Delta _{{\mathop{\rm Im}\nolimits} }}$ are within the range $\pm 0.001$. Further increase of the integer $N$ significantly improves the accuracy. This can be seen from Fig. 7 showing that the differences ${\Delta _{{\mathop{\rm Re}\nolimits} }}$ (blue curve) and ${\Delta _{{\mathop{\rm Im}\nolimits} }}$ (red curve) computed at $N = 300$, $h = c = 0.00166389$ remain within the narrow range $ \pm 0.00003$.

\begin{figure}[ht]
\begin{center}
\includegraphics[width=24pc]{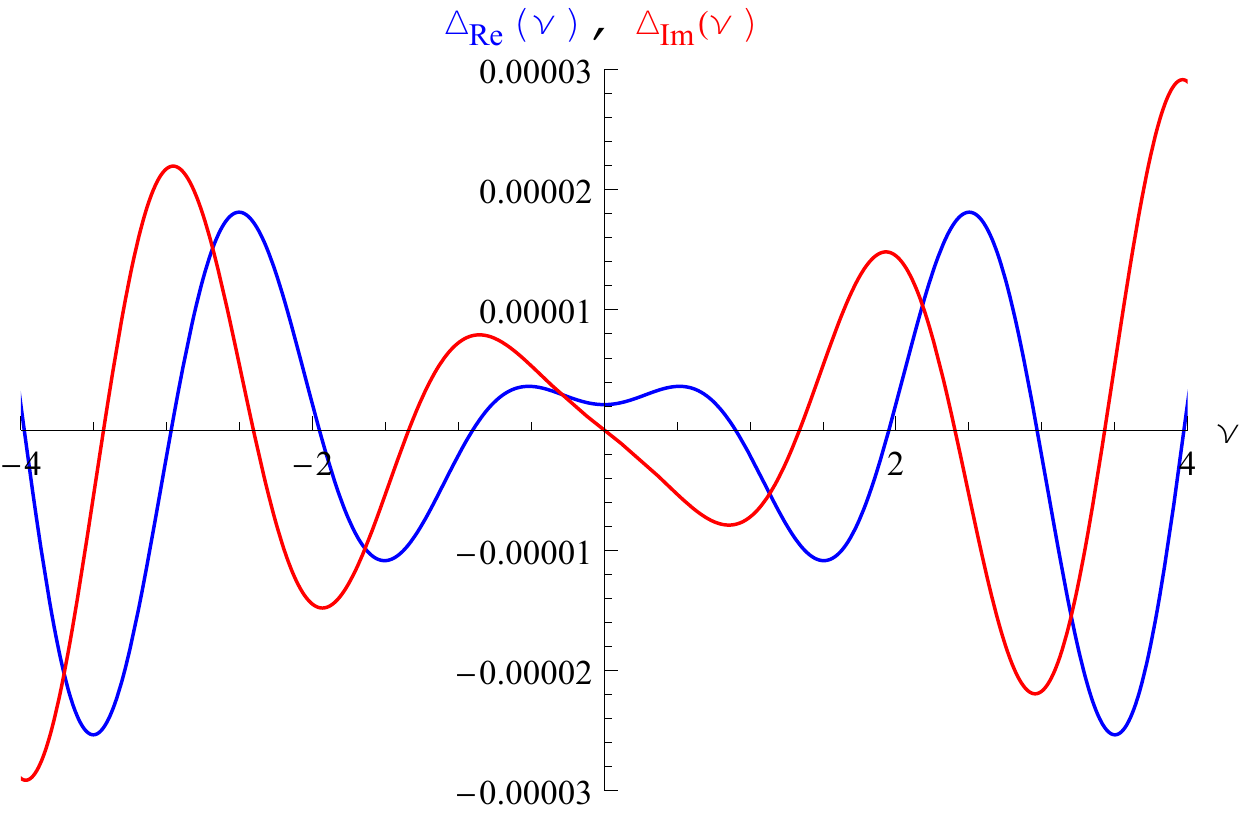}\hspace{2pc}%
\begin{minipage}[b]{28pc}
\vspace{0.3cm}
{\sffamily {\bf{Fig. 7.}} The differences ${\Delta _{{\mathop{\rm Re}\nolimits} }}$ (blue curve) and ${\Delta _{{\mathop{\rm Im}\nolimits} }}$ (red curve) computed at $N = 300$, $h = c = 0.00166389$, respectively.}
\end{minipage}
\end{center}
\end{figure}

\section{Conclusion}

We present a new approach for numerical computation of the Fourier integrals based on a sampling with the Gaussian function of kind $h\,{e^{ - {{\left( {t/c} \right)}^2}}}/\left( {{c}\sqrt \pi  } \right)$. It is shown that the Fourier transform can be expressed as a weighted sum of the complex error functions. Applying a remarkable property of the complex error function, the weighted sum of the complex error functions can be significantly simplified as a damping harmonic series. Unlike the conventional discrete Fourier transform, this methodology results in a non-periodic wavelet approximation. Therefore, the proposed approach may be practically convenient and advantageous in algorithmic implementation.

\section*{Acknowledgments}

This work is supported by National Research Council Canada, Thoth Technology Inc. and York University. The authors wish to thank to Prof. Ian McDade and Dr. Brian Solheim for discussions and constructive suggestions.


\end{document}